\input amstex
\input amsppt.sty
\magnification=\magstep1
\hsize=30truecc
\vsize=22.2truecm
\baselineskip=16truept
\TagsOnRight \pageno=1 \nologo

\def\N{\Bbb N}

\def\l{\left}
\def\r{\right}
\def\bg{\bigg}
\def\({\bg(}
\def\[{\bg\lfloor}
\def\){\bg)}
\def\]{\bg\rfloor}
\def\t{\text}
\def\f{\frac}

\def\bi{\binom}
\def\eq{\equiv}

\def\ls{\leqslant}

\def\mo{\roman{mod}}

\def\da{\delta}

\def\Proof{\noindent{\it Proof}}

\hbox {Proc. Amer. Math. Soc. 138(2010), no.\,1, 37--46.}
\bigskip
\topmatter
\title Some congruences for the second-order Catalan numbers\endtitle
\author Li-Lu Zhao, Hao Pan and Zhi-Wei Sun\endauthor
\leftheadtext{Li-Lu Zhao, Hao Pan and Zhi-Wei Sun}
\affil Department of Mathematics, Nanjing University\\
 Nanjing 210093, People's Republic of China
\\{\tt zhaolilu\@gmail.com, haopan79\@yahoo.com.cn,
zwsun\@nju.edu.cn}
\endaffil
\abstract Let $p$ be any odd prime. We mainly show that
$$\sum_{k=1}^{p-1}\f{2^k}k\bi{3k}k\eq0\ (\mo\ p)$$
and
$$\sum_{k=1}^{p-1}2^{k-1}C_k^{(2)}\eq(-1)^{(p-1)/2}-1\ (\mo\ p),$$
where $C_k^{(2)}=\bi{3k}k/(2k+1)$ is the $k$th Catalan number of order 2.
\endabstract
\thanks 2000 {\it Mathematics Subject Classification}.\,Primary 11B65;
Secondary 05A10,\,11A07.
\newline\indent The third author is the corresponding author.
He was supported by the National Natural Science
Foundation (grant 10871087) and the Overseas Cooperation Fund of China.
\endthanks
\endtopmatter
\document

\heading{1. Introduction}\endheading

 The well-known Catalan numbers are those integers
$$C_n=\f 1{n+1}\bi{2n}n=\bi{2n}n-\bi{2n}{n-1} \ \ (n=0,1,2,\ldots).$$
(As usual we regard $\bi x{-k}$ as 0 for $k=1,2,\ldots$.)
There are many combinatorial interpretations for these important
numbers (see, e.g., [St, pp. 219-229]). With the help of a
sophisticated binomial identity, H. Pan and Z. W. Sun [PS]
obtained some congruences on sums of Catalan numbers; in
particular, by [PS, (1.16) and (1.8)], for any prime $p>3$ we have
$$\sum_{k=0}^{p-1}C_k\eq\f{3(\f p3)-1}2\ (\mo\ p)
\ \ \t{and}\ \ \sum_{k=1}^{p-1}\f{C_k}k\eq\f32\l(1-\l(\f p3\r)\r)\
(\mo\ p),\tag1.0$$ where the Legendre symbol $(\f
a3)\in\{0,\pm1\}$ satisfies the congruence $a\eq(\f a3)\ (\mo\
3)$. Recently Z. W. Sun and R. Tauraso [ST1, ST2] obtained some further congruences concerning sums
involving Catalan numbers.

For $m,n\in \N=\{0,1,2,\ldots\}$, we define
$$C_n^{(m)}=\f {1}{mn+1}\bi{mn+n}n=\bi{mn+n}n-m\bi{mn+n}{n-1}$$
and call it the $n$th Catalan number of order $m$. Clearly
$$C_n^{(1)}=C_n\quad\ \t{and}\quad\ C_n^{(2)}=\f{1}{2n+1}\bi{3n}n.$$
In contrast with (1.0), we have the following result involving the
second-order Catalan numbers.

\proclaim{Theorem 1.1} Let $p$ be an odd prime. Then
$$\sum_{k=1}^{p-1}2^{k}C_k^{(2)}\eq 2\l((-1)^{(p-1)/2}-1\r) \ (\mo\ p)\tag1.1$$
and
$$\sum_{k=1}^{p-1} \f  {2^{k}C_k^{(2)}}{k} \eq 4\l(1-(-1)^{(p-1)/2}\r)\ (\mo\ p).\tag1.2$$
\endproclaim

Actually Theorem 1.1 follows from our next two theorems.

\proclaim{Theorem 1.2} Let $p>5$ be a prime. Then
$$\align\sum_{k=0}^{p-1}2^k\bi{3k}k\eq&\f{6(-1)^{(p-1)/2}-1}5\ (\mo\ p),\tag1.3
\\ \sum_{k=0}^{p-1}2^k\bi{3k+1}k\eq&\f{4(-1)^{(p-1)/2}+1}5\ (\mo\ p),\tag1.4
\endalign$$
\endproclaim

\proclaim{Theorem 1.3} For any prime $p$ we have
$$\sum_{k=1}^{p-1}\f{2^k}k\bi{3k}k\eq0\ (\mo\ p).\tag1.5$$
\endproclaim

 For any odd prime $p$ we can also prove the following congruences:
 $$\align 5\sum_{k=1}^{p-1}2^k\bi{3k+2}k\eq&(-1)^{(p-1)/2}-1\ (\mo\ p),
 \\\sum_{k=1}^{p-1}\f{2^{k-1}}k\bi{3k+1}k\eq&(-1)^{(p-1)/2}-1\ (\mo\ p),
 \\\sum_{k=1}^{p-1}\f{2^{k-1}}k\bi{3k+2}k\eq&\f
 32\l((-1)^{(p-1)/2}-1\r)\ (\mo\ p).
 \endalign$$
We omit their proofs which are similar to those of Theorems 1.2-1.3.

With the help of Theorems 1.2 and 1.3, we can easily deduce Theorem
1.1.

\noindent{\it Proof of Theorem 1.1 via Theorems 1.2 and 1.3}.
Clearly (1.1) and (1.2) hold for $p=3,5$. Assume $p>5$. By (1.3) and
(1.4),
$$\aligned
\sum_{k=0}^{p-1} \f {2^k}{2k+1} \bi {3k}k=&
3\sum_{k=0}^{p-1}2^k\bi{3k}k-2\sum_{k=0}^{p-1} 2^k \bi {3k+1}k
\\\eq&2(-1)^{(p-1)/2}-1\ (\mo\ p).
\endaligned
$$
This proves (1.1). For (1.2) it suffices to note that
$$\sum_{k=1}^{p-1} \f {2^k}{k(2k+1)} \bi {3k}k
=\sum_{k=1}^{p-1} \f{2^k}{k} \bi {3k}k-2\sum_{k=1}^{p-1} \f {2^k}{2k+1} \bi {3k}k.$$
This concludes the proof. \qed

We are going to provide two lemmas in the next section.
Theorems 1.2 and 1.3 will be proved in Sections 3 and 4 respectively.

 \heading{2. Some Lemmas}\endheading

\proclaim{Lemma 2.1} For $m,n\in\N$ we have
$$\aligned& 2^n \sum_{k=0}^{\lfloor m/3\rfloor} (-2)^k\bi {n}{m-3k} \bi {3k-m+n}k
\\=&(-1)^m \sum_{j=0}^{n}\bi nj\sum_{k=0}^{m}(-2)^k\bi n{m-k}\bi {2j}k. \endaligned\tag2.1$$
\endproclaim
\Proof. Let $P(x)=(2+2x-4x^3)^n$, and denote by $[x^k]P(x)$ the
coefficient of $x^k$ in the expansion of $P(x)$. Then
$$\aligned 2^{-n}[x^m]P(x)=&[x^m]((1+x)-2x^3)^n
\\=&\sum_{k=0}^{\lfloor m/3\rfloor} \bi nk(-2)^k[x^{m-3k}](1+x)^{n-k}\\
=&\sum_{k=0}^{\lfloor m/3\rfloor} (-2)^k\bi
n{k}\bi{n-k}{m-3k}\\
=&\sum_{k=0}^{\lfloor m/3\rfloor} (-2)^k
\bi{n}{m-3k}\bi{3k-m+n}{k}.
\endaligned$$
Since
$$P(x)=(1-x)^n((2x+1)^2+1)^n=\sum_{j=0}^{n}
\bi nj(1-x)^n(2x+1)^{2j},$$ we also have
$$[x^m]P(x)=\sum_{j=0}^{n}\bi nj\sum_{k=0}^{m}
2^k\bi {2j}k(-1)^{m-k}\bi{n}{m-k}.$$ Therefore (2.1) is valid.
\qed

For any prime $p$, if $n,k\in\N$ and $s,t\in\{0,1,\ldots,p-1\}$ then we have
the following well-known Lucas congruence (cf. [Gr] or [HS]):
$\bi {pn+s}{pk+t}\eq\bi nk\bi st\ (\mo\ p)$.
This will be used in the proof of the following lemma.

\proclaim{Lemma 2.2} Let $p>5$ be a prime. Then we have
$$\sum_{s=0}^{p-1}(-1)^s\sum_{t=0}^{p-1} 2^{t}\bi{2s}{t}\eq\f{3(-1)^{(p-1)/2}+2}{5}\ (\mo\ p)\tag2.2$$
and
$$\sum_{s=0}^{p-1}(-1)^s\sum_{t=0}^{p-1} 2^{t}\bi{2s}{p+t}\eq\f{3}{10}\l(1-(-1)^{(p-1)/2}\r)\ (\mo\ p).\tag2.3$$
\endproclaim
\Proof. Observe that
$$\align&\sum_{s=0}^{p-1}(-1)^s\sum_{t=0}^{p-1}
2^{t}\bi{2s}{t}
\\=&\sum_{s=0}^{(p-1)/2}(-1)^s\sum_{t=0}^{2s}
2^{t}\bi{2s}{t}+\sum_{s=(p+1)/2}^{p-1}(-1)^s\sum_{t=0}^{p-1}
2^{t}\bi{2s}{t}
\\=&\sum_{s=0}^{(p-1)/2}(-1)^s3^{2s}+\sum_{s=(p+1)/2}^{p-1}(-1)^s\sum_{t=0}^{p-1}
2^{t}\bi{2s}{t}
\\=&\sum_{s=0}^{(p-1)/2}(-1)^s3^{2s}+\sum_{s=(p+1)/2}^{p-1}(-1)^s\(\sum_{t=0}^{2s}
2^{t}\bi{2s}{t}-\sum_{t=p}^{2s} 2^{t}\bi{2s}{t}\)
\\=&\sum_{s=0}^{p-1}(-1)^s3^{2s}-\sum_{s=(p+1)/2}^{p-1}(-1)^s\sum_{t=p}^{2s} 2^{t}\bi{2s}{t}
\\=&\sum_{s=0}^{p-1}(-9)^s-\sum_{s=(p+1)/2}^{p-1}(-1)^s\sum_{r=0}^{2s-p} 2^{p+r}\bi{2s}{r+p}.
\endalign$$
For $s=(p+1)/2,\ldots,p-1$, by Lucas' congruence  we have
$$\sum_{r=0}^{2s-p}2^r\bi{p+(2s-p)}{p+r}\eq\sum_{r=0}^{2s-p}2^r\bi{2s-p}r=3^{2s-p}\ (\mo\ p).$$
Thus, with the help of Fermat's little theorem, we get
$$\align\sum_{s=0}^{p-1}(-1)^s\sum_{t=0}^{p-1}
2^{t}\bi{2s}{t}
\eq&\f{1-(-9)^p}{10}-\sum_{s=(p+1)/2}^{p-1}(-1)^s\f 23\cdot 9^s
\\\eq&1-\f{2}{3}(-9)^{\f{p+1}{2}}\f{(1-(-9)^{(p-1)/2})}{10}
\\\eq&\f{3(-1)^{(p-1)/2}+2}{5}\ (\mo\ p).
\endalign$$
This proves (2.2).

In view of Lucas' congruence and Fermat's little theorem, we also have
$$\align&\sum_{s=0}^{p-1}(-1)^s\sum_{t=0}^{p-1}2^t\bi{2s}{p+t}
\\\eq&\sum_{s=(p+1)/2}^{p-1}(-1)^s\sum_{t=0}^{p-1}2^t\bi{2s-p}t=\sum_{s=(p+1)/2}^{p-1}(-1)^s3^{2s-p}
\\=&3^{-p}(-9)^{(p+1)/2}\f{1-(-9)^{(p-1)/2}}{10}=(-1)^{(p+1)/2}\f 3{10}\l(1+(-1)^{(p+1)/2}3^{p-1}\r)
\\\eq&\f{3}{10}\l(1-(-1)^{(p-1)/2}\r)\ (\mo\ p).
\endalign$$
So (2.3) is also valid. We are done. \qed

\heading{3. Proof of Theorem 1.2}\endheading

In order to prove Theorem 1.2, we first present an auxiliary result.

\proclaim{Theorem 3.1} Let $p>5$ be a prime, and let
$d,\da\in\{0,1\}$. Then
$$\aligned&\f{(-1)^{d+\da}}{2^{\da}}\sum_{\da p-d\ls 3k\ls \da p+p-1-d} 2^k \bi {3k+d}k
\\\eq&\f{4-\da}{10}+\f{(3\da-2)(5d-3)}{10}(-1)^{(p-1)/2}\ (\mo\ p).
\endaligned\tag3.1$$
\endproclaim
\Proof. Applying (2.1) with $n=p-1$
and $m=\da p+p-1-d$, we get
$$\align&2^{p-1} \sum_{k=0}^{\lfloor (\da p+p-1-d)/3\rfloor} (-2)^k\bi {p-1}{\da p+p-1-d-3k} \bi {3k+d-\da p}k
\\=&(-1)^{\da p+p-1-d}\sum_{j=0}^{p-1} \bi {p-1}j\sum_{k=0}^{\da p+p-1-d}
(-2)^k\bi {p-1}{\da p+p-1-d-k}\bi{2j}k.
\endalign$$
Observe that
$$\align&\sum_{k=0}^{\lfloor (\da p+p-1-d)/3\rfloor} (-2)^k\bi {p-1}{\da p+p-1-d-3k} \bi {3k+d-\da p}k
\\=&\sum_{\da p-d\ls 3k\ls \da p+p-1-d} (-2)^k\bi {p-1}{p+\da p-1-d-3k} \bi {3k+d-\da p}k
\\\eq& \sum_{\da p-d\ls 3k\ls \da p+p-1-d} (-2)^k(-1)^{\da p+p-1-d-3k} \bi {3k+d}k
\\\eq&(-1)^{d+\da}\sum_{\da p-d\ls 3k\ls \da p+p-1-d} 2^k \bi {3k+d}k\ (\mo\ p)
\endalign$$
and
$$\align&(-1)^{\da p+p-1-d}\sum_{j=0}^{p-1} \bi {p-1}j\sum_{k=0}^{\da p+p-1-d}
(-2)^k\bi {p-1}{\da p+p-1-d-k}\bi{2j}k
\\\eq&\sum_{j=0}^{p-1}(-1)^j\sum_{\da p-d\ls k<\da p+p-d} 2^k\bi{2j}k
=\sum_{j=0}^{p-1}(-1)^j\sum_{t=0}^{p-1} 2^{\da p-d+t}\bi{2j}{\da p-d+t}
\\\eq&2^{\da-d}\sum_{s=0}^{p-1}(-1)^s\sum_{t=0}^{p-1} 2^{t}\bi{2s}{\da p-d+t}\ (\mo\ p).
\endalign$$
Therefore
$$\align&\sum_{\da p-d\ls 3k\ls \da p+p-1-d} 2^k \bi {3k+d}k
\\\eq&(-2)^{\da-d}\sum_{s=0}^{p-1}(-1)^s\sum_{t=0}^{p-1} 2^{t}\bi{2s}{\da p-d+t}\ (\mo\ p).
\endalign$$

 Recall that $d\in\{0,1\}$. We have
$$\align&\sum_{s=0}^{p-1}(-1)^s\sum_{t=0}^{p-1} 2^{t}\bi{2s}{\da p-d+t}
\\=&\sum_{s=0}^{p-1}(-1)^s\sum_{t=-d}^{p-1-d} 2^{d+t}\bi{2s}{\da p+t}
\\=&\sum_{s=0}^{p-1}(-1)^s\(\sum_{t=0}^{p-1} 2^{d+t}\bi{2s}{\da p+t}+d\(\bi{2s}{\da p-1}-2^{p}\bi{2s}{\da p+p-1}\)\)
\\=&2^d\sum_{s=0}^{p-1}\sum_{t=0}^{p-1}(-1)^s2^t\bi{2s}{\da p+t}
+d\sum_{s=0}^{p-1}(-1)^s\(\bi{2s}{\da p-1}-2^p\bi{2s}{\da p+p-1}\)
\endalign$$
and hence
$$\align&(-1)^{d+\da}\sum_{\da p-d\ls 3k\ls \da p+p-1-d} 2^k \bi {3k+d}k
-2^\da\sum_{s=0}^{p-1}\sum_{t=0}^{p-1}(-1)^s2^t\bi{2s}{\da p+t}
\\\eq&d2^{\da-d}\sum_{s=0}^{p-1}(-1)^s\(\bi{2s}{\da p-1}-2\bi{2s}{\da p+p-1}\)
\\\eq&d2^{\da-1}(3\da-2)\sum_{s=0}^{p-1}(-1)^s\bi{2s}{p-1}\eq d(3\da-2)2^{\da-1}(-1)^{(p-1)/2}\ (\mo\ p).
\endalign$$
Since
$$\sum_{s=0}^{p-1}(-1)^s\sum_{t=0}^{p-1}2^t\bi{2s}{\da p+t}\eq\f{4-\da}{10}+\f3{10}(2-3\da)(-1)^{(p-1)/2}\ (\mo\ p)$$
by Lemma 2.2, we finally get
$$\align&\f{(-1)^{d+\da}}{2^{\da}}\sum_{\da p-d\ls 3k\ls \da p+p-1-d} 2^k \bi {3k+d}k
\\\eq&\f{4-\da}{10}+\f3{10}(2-3\da)(-1)^{(p-1)/2}+\f d2(3\da-2)(-1)^{(p-1)/2}
\\\eq&\f{4-\da}{10}+\f{(3\da-2)(5d-3)}{10}(-1)^{(p-1)/2}\ (\mo\ p).
\endalign$$
This proves (3.1). \qed

\medskip
\noindent{\it Proof of Theorem 1.2}. Let $d\in\{0,1\}$. If $(2p-d)/3\ls k\ls p-1$, then
$2k+d+1\ls 2k+2\ls 2p\ls 3k+d$ and hence
$$\bi{3k+d}k=\f{(3k+d)\cdots(2k+d+1)}{k!}\eq0\ (\mo\ p).$$
Therefore
$$\sum_{2p-d\ls 3k\ls 3p-3}2^k\bi{3k+d}k\eq0\ (\mo\ p).$$
With the help of Theorem 3.1, we have
$$\align \sum_{k=0}^{p-1}2^k\bi{3k+d}k\eq &\sum_{-d\ls 3k\ls 2p-1-d}2^k\bi{3k+d}k
\\\eq&\sum_{\da=0}^1\ \sum_{\da p-d\ls 3k\ls \da p+p-1-d}2^k\bi{3k+d}k
\\\eq&\sum_{\da=0}^1(-1)^d(-2)^{\da}\(\f{4-\da}{10}+\f{(3\da-2)(5d-3)}{10}(-1)^{(p-1)/2}\)
\\\eq&\f{(-1)^{d-1}}5\l(1+(10d-6)(-1)^{(p-1)/2}\r)\ (\mo\ p).
\endalign$$
This yields (1.3) and (1.4). We are done. \qed

\heading{4. Proof of Theorem 1.3}\endheading

\medskip
\noindent{\it Proof of Theorem 1.3}. Obviously (1.5) holds for $p=2,3$. Below we assume $p>3$.

Let $\da\in\{0,1\}$.
Applying (2.1) with $m=p+\da p$ and $n=p$  we get
$$\aligned& 2^{p} \sum_{k=0}^{p} (-2)^k\bi {p}{p+\da p-3k} \bi {3k-\da p}k
\\=&(-1)^{\da+1}\sum_{j=0}^{p} \bi pj\sum_{k=0}^{p+\da p}(-2)^k\bi {p}{p+\da p-k}\bi{2j}k.
\endaligned\tag4.1$$

 Observe that
$$\aligned& \sum_{k=0}^{p} (-2)^k\bi {p}{p+\da p-3k} \bi {3k-\da p}k
\\ =&\sum_{\da p\ls 3k\ls p+\da p-1} (-2)^k\bi {p}{3k-\da p} \bi {3k-\da p}k
\\=&1-\da+\sum_{\da p<3k<p+\da p}(-2)^k\bi p{3k-\da p}\bi{3k-\da p}k.
\endaligned$$
For $j=1,\ldots,p-1$ clearly
$$\bi pj=\f pj\bi{p-1}{j-1}\eq p\f{(-1)^{j-1}}j\ (\mo\ p^2).$$
Thus
$$\align &\sum_{\da p<3k<p+\da p}(-2)^k\bi p{3k-\da p}\bi{3k-\da p}k
\\\eq&\sum_{\da p<3k<p+\da p}(-2)^k p\f{(-1)^{3k-\da p-1}}{3k-\da p}\bi{3k-\da p}k
\\\eq&(-1)^{\da+1}\sum_{\da p<3k<p+\da p}(-2)^k p\f{(-1)^{k}}{3k}\bi{(3k-\da p)+\da p}k
\\&\qquad\qquad\qquad\qquad\qquad\qquad\qquad\qquad\ (\t{by Lucas' congruence})
\\\eq&(-1)^{\da+1}\ \f p3\sum_{\da p<3k<p+\da p}\f{2^k}k\bi{3k}k\ (\mo\ p^2).
\endalign$$

 Notice that
$$\align &\sum_{j=0}^{p}\bi {p}j \sum_{k=0}^{p+\da p}(-2)^k\bi {p}{p+\da p-k}\bi
{2j}k\\
=&\sum_{\da p\ls 2j\ls 2p}\bi {p}j \sum_{k=\da p}^{p+\da p}(-2)^k\bi {p}{k-\da p}\bi
{2j}k\\
=&\sum_{\da p<2j<2p}\bi {p}j \sum_{k=\da p}^{p+\da p}(-2)^k\bi
{p}{k-\da p}\bi {2j}k
\\&+\sum_{2j\in\{\da p,2p\}}\bi pj\sum_{k=\da p}^{p+\da p}(-2)^k\bi {p}{k-\da p}\bi {2j}k.
\endalign$$
Clearly
$$\align&\sum_{\da p<2j<2p}\bi {p}j \sum_{k=\da p}^{p+\da p}(-2)^k\bi
{p}{k-\da p}\bi {2j}k
\\\eq&\sum_{\da p<2j<2p}\bi pj\((-2)^{\da p}\bi p0\bi{2j}{\da p}+(-2)^{p+\da p}\bi pp\bi{2j}{p+\da p}\)
\\\eq&\sum_{\da p<2j<2p}\bi pj(-2)^{\da p}\bi{2j-\da p}{0}
\\&+(1-\da)\sum_{p<2j<2p}\bi pj(-2)^{p+\da p}\bi{2j-p}{p-p}\ (\t{by Lucas' congruence})
\\\eq&(-2)^{\da}2^{1-\da}(2^{p-1}-1)+(1-\da)(-2)^{1+\da}(2^{p-1}-1)
\\\eq&(-1)^{\da}\da(2^p-2)=-\da(2^p-2)\ (\mo\ p^2).
\endalign$$
(Note that $\da\in\{0,1\}$ and $2\sum_{p/2<j<p}\bi pj=\sum_{j=1}^{p-1}\bi pj=2^p-2$.)
Also,
$$\sum_{2j=\da p}\bi pj\sum_{k=\da p}^{p+\da p}(-2)^k\bi p{k-\da p}\bi{2j}k
=(1-\da)\sum_{k=0}^p(-2)^k\bi pk\bi0k=1-\da$$
and
$$\align&\sum_{2j=2p}\bi pj\sum_{k=\da p}^{p+\da p}(-2)^k\bi p{k-\da p}\bi{2j}k
\\\eq&\sum_{k\in\{\da p,p+\da p\}}(-2)^k\bi p{k-\da p}\bi{2p}k
\\\eq&(-2)^{\da p}\bi 2\da+(-2)^{p+\da p}\bi 2{1+\da}=4^{\da p}-2^{p+1}\ (\mo\ p^2).
\endalign$$
(Recall that $\f12\bi {2p}p=\bi{2p-1}{p-1}\eq1\ (\mo\ p^3)$ by the Wolstenholme congruence (cf. [Gr] or [HT]).)

Combining the above with (4.1), we have
$$\align&2^p\(1-\da+(-1)^{\da+1}\ \f p3\sum_{\da p<3k<p+\da p}\f{2^k}k\bi{3k}k\)
\\\eq&(-1)^{\da+1}\l(\da(2-2^p)+1-\da+4^{\da p}-2^{p+1}\r)\ (\mo\ p^2).
\endalign$$
Setting $\da=0$ and $\da=1$ respectively, we obtain
$$2^p-2^p\f p3\sum_{0<3k<p}\f{2^k}k\bi{3k}k\eq 2^{p+1}-2\ (\mo\ p^2)$$
and
$$2^p\f p3\sum_{p<3k<2p}\f{2^k}k\bi{3k}k\eq 2-2^p+4^p-2^{p+1}\ (\mo\ p^2).$$
It follows that
$$\f 23p\sum_{0<3k<2p}\f{2^k}k\bi{3k}k\eq 4^p-4\cdot 2^p+4=(2^p-2)^2\eq0\ (\mo\ p^2).$$
If $2p\ls 3k<3p$, then
$$\bi{3k}k=\f{3k\cdots(2k+1)}{k!}\eq0\ (\mo\ p).$$
Therefore
$$\sum_{k=1}^{p-1}\f{2^k}k\bi{3k}k=\sum_{0<3k<2p}\f{2^k}k\bi{3k}k+\sum_{2p\ls 3k<3p}\f{2^k}k\bi{3k}k\eq0\ (\mo\ p).$$
This completes the proof of Theorem 1.3. \qed

\enddocument